\documentclass{ifacconf}

\usepackage[caption=false]{subfig}
\usepackage[dvips]{graphicx}
\graphicspath{{figures/}}
\usepackage{amsmath, amssymb}
\usepackage{tikz} \usetikzlibrary{calc}
\usepackage{psfrag}
\usepackage{pgfplots}
\usepackage{algorithm}     % <-- should be replaced with 'algorithms'
\usepackage{algpseudocode} % <-- should be replaced with 'algorithms'
\usepackage{arydshln}

\DeclareMathOperator*{\minimize}{minimize}

\DeclareMathOperator*{\subject}{subj.\ to}

\DeclareMathOperator*{\diag}{diag}

\DeclareMathOperator*{\find}{find}
\DeclareMathOperator*{\adj}{adj}
\DeclareMathOperator*{\parent}{par}	
\DeclareMathOperator*{\IQC}{IQC}
\DeclareMathOperator*{\children}{ch}
\DeclareMathOperator*{\degree}{deg}

\definecolor{red}{rgb}{1,0,0}

\begin{document}

\begin{frontmatter}

\title{Distributed Robustness Analysis of Interconnected Uncertain Systems Using Chordal Decomposition\thanksref{footnoteinfo}} % Title, preferably not more than 10 words.
\vspace*{-8pt}
\thanks[footnoteinfo]{This work has been supported by the Swedish Department of Education within the ELLIIT project.}
\author[First]{Sina Khoshfetrat Pakazad} 
\author[First]{Anders Hansson}
\author[Second]{Martin S. Andersen} 
\author[Third]{Anders Rantzer}
\vspace*{-8pt}
\address[First]{Division of Automatic Control, Department of Electrical Engineering, Link\"oping University, Sweden. Email: \{sina.kh.pa, hansson\}@isy.liu.se.} 
\address[Second]{Department of Applied Mathematics and Computer Science, Technical University of Denmark, Denmark. (e-mail: mskan@dtu.dk)}  
\address[Third]{Department of Automatic Control, Lund University, Sweden. (e-mail: anders.rantzer@control.lth.se)} 
\vspace*{-8pt}

%\begin{keyword}                           % Five to ten keywords,  
%Cicero; Catiline; orations.               % chosen from the IFAC 
%\end{keyword}                             % keyword list or with the 
                                          % help of the Automatica 
                                          % keyword wizard

\begin{abstract}                          % Abstract of not more than 250 words.
Large-scale interconnected uncertain systems commonly have large state and uncertainty dimensions. Aside from the heavy computational cost of solving centralized robust stability analysis techniques, privacy requirements in the network can also introduce further issues. In this paper, we utilize IQC analysis for analyzing large-scale interconnected uncertain systems and we evade these issues by describing a decomposition scheme that is based on the interconnection structure of the system. This scheme is based on the so-called chordal decomposition and does not add any conservativeness to the analysis approach. The decomposed problem can be solved using distributed computational algorithms without the need for a centralized computational unit. We further discuss the merits of the proposed analysis approach using a numerical experiment.
\end{abstract}

\end{frontmatter}
%================================================================
\section{Introduction}
%================================================================
\vspace*{-8pt}
Stability analysis of uncertain systems investigates whether a system is stable with respect to all the admissible values of uncertain parameters or not. In case the uncertain system is stable with respect to all admissible uncertainties, it is said that the system is robustly stable, \cite{robustandoptimal}. There are different approaches for analyzing robust stability of uncertain systems and IQC analysis is one of the more general ones, \cite{ulfiqc,rantzer}. Analyzing robustness using IQC analysis usually requires solving a semi-infinite frequency-dependent linear matrix inequality (LMI), and there are mainly two approaches for solving this LMI which either rely on the use of the  KYP lemma or on frequency gridding. These approaches in turn require solving a finite number of LMIs, the dimension and number of variables of which usually grow rapidly with the dimension of the system, i.e., number of states and/or number of input-outputs, and the size of the uncertainty block in the system. This makes analyzing high-dimensional systems computationally expensive and even prohibitive. Such computational issues can be alleviated if we exploit the inherent structure in the resulting LMIs, see e.g., \cite{ChapterBook1,Anders,Parrilo,Kao2,and+dah+van10,WHJ:09,Kao1,and:13}. However the computational burden of solving the analysis problem can still remain prohibitive, specially when we intend to analyze interconnections of uncertain subsystems, with large number of subsystems. Furthermore, the mentioned methods become unviable if we have local privacy requirements in the network, i.e., the subsystems in the network are not willing to share their system model with a central computational unit and/or with all other subsystems. This prevents us from being able to form and solve the problem in a centralized manner. Both these issues can be evaded through decomposition of the problem and utilizing distributed algorithms for solving it, \cite{UlfLetter,Kao,Ulf,Langbort,P5,Vinnicombe,kim:12}. The authors in \cite{UlfLetter} utilize IQC analysis to analyze robustness of interconnected uncertain systems when the interconnection among the subsystems is described as $\Gamma = \bar \Gamma \otimes I$, where $\bar \Gamma$ is the adjacency matrix of the network. They then propose a decomposition for the analysis problem. This decomposition, however, only applies when $\bar \Gamma$ is normal. The authors in \cite{Kao,Ulf} also consider analysis of interconnected systems, though, when the uncertainty lies in the interconnections. They also provide decomposition schemes for the analysis problem, when the interconnections are described as $\Gamma = \bar \Gamma \otimes I$ with $\bar \Gamma$ being normal. However, these decomposition schemes are only valid when the spectrum of $\bar \Gamma$ can either be characterized using quadratic, \cite{Kao}, or polyhedral, \cite{Ulf}, constraints. Similar problems are also considered in \cite{Langbort,P5}, for more general interconnections. However, the computational complexity or details on how to solve the decomposed problem are not discussed in these papers. Decomposition and distributed solutions for analyzing interconnections of uncertain homogeneous, i.e., when all the subsystems have the same model, and heterogeneous SISO subsystems are considered in \cite{kim:12,Vinnicombe}, respectively.\newline
The methods in the above mentioned papers, provide viable approaches for decomposing the analysis of interconnected uncertain systems. However, they all rely on strong assumptions on the interconnections, e.g., normality of the adjacency matrix of the network and polyhedral or quadratic characterization of spectrum of this matrix, and/or on the underlying subsystems in the network, e.g., their input-output dimensions and homogeneity of subsystems. In this paper, we do not make strong assumptions regarding the interconnections or subsystems and we merely assume that the interconnections among subsystems are sparse, i.e., each subsystem is connected to only a few other subsystems in the network. This is usually a characteristic of large-scale interconnected uncertain systems. In \cite{and:13}, we used IQC analysis for analyzing large-scale interconnected uncertain systems and we proposed an equivalent reformulation of the analysis problem which enabled us to reflect the inherent sparsity in interconnections among subsystems in the resulting semi-infinite LMI. We then showed that the computational cost of analyzing large-scale interconnected uncertain systems using frequency-gridding, in a centralized manner, can be reduced significantly if we exploit the structure imposed by inherent sparsity in the interconnection. However, as was mentioned earlier, centralized solutions can be unviable for analyzing the system. So in \cite{and+han:12}, using the same formulation, we put forth the preliminary ideas on how to decompose the analysis problem using its sparsity pattern. \newline
In this paper, we also utilize the formulation proposed in \cite{and:13,and+han:12}, and in addition we
\begin{itemize}
\item provide a more efficient way to handle interconnected uncertain systems with uncertain interconnections which also gives us further insight on how to analyze interconnected uncertain systems without uncertain interconnectionsš
\item describe the decomposition process of analyzing interconnected uncertain systems with certain and uncertain interconnections;
\item discuss how the decomposed problem can be solved efficiently in a distributed manner, and how this can enable us to provide partial local privacy.  
\end{itemize}
\vspace*{-7pt}
\subsection*{Notation}
\vspace*{-8pt}
The set of real scalars and $m \times n$ real matrices are denoted by $\mathbb R$ and $\mathbb R^{m\times n}$, respectively. We denote transpose  of a matrix $G$ by $G^T$ and its conjugate transpose by $G^{\ast}$. The set $\mathbf S^n$ denotes $n \times n$ Hermitian matrices, and $I_n$ denotes the $n \times n$ identity matrix. Given a set $N$ with $|N|$ we denote the number of elements in the set, and we denote the set-theoretic difference between two sets $N_1$ and $N_2$ by $N_1\setminus N_2$. We use superscripts for indexing different matrices, and we use subscripts to refer to different elements in the matrix, i.e., by $G^{k}_{ij}$ denotes the element on the $i$th row and $j$th column of the matrix $G^k$. Similarly, $v^k_i$ is the $i$th component of the vector $v^k$. Given matrices $G^k$ for $k = 1, \dots, N$, by $\diag(G^1, \dots, G^N)$ we denote a block-diagonal matrix with diagonal blocks specified by the given matrices. Likewise, given vectors $v^k$ for $k= 1, \dots, N$, the column vector $(v^1, \dots, v^N)$ is all of the given vectors stacked. By $\mathcal L_2^n$ we denote the set of $n$-dimensional square integrable signals, and $\mathcal{RH}_{\infty}^{m \times n}$ represents the set of real, rational $m \times n$ transfer function matrices with no poles in the closed right half plane. 
\vspace*{-8pt}
%=================================================================
\section{Robust Stability Analysis of Uncertain Systems Using IQCs}\label{sec:IQC}
%=================================================================
\vspace*{-8pt}
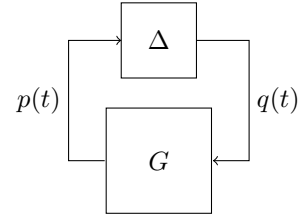
\begin{figure}
\begin{center}
\begin{tikzpicture}[scale=0.4, inner sep=0pt]

  \draw (0,0) node [draw, fill=none, minimum size=1.4cm]
  {$G$};
  \draw (0,4)  node [draw, fill=none, minimum size=1.0cm]
  {$\Delta$};

  \draw[->] (-1.8,0) -- (-3,0) -- node[left]
  {$p(t)\phantom{.}$} (-3,4) -- (-1.25,4);
  \draw[->] (1.25,4) -- (3,4) -- node[right]
  {$\phantom{.}q(t)$} (3,0) -- (1.8,0);

\end{tikzpicture}  
\vspace*{-7pt}
\caption{\small Uncertain system with system transfer function matrix $G$ and uncertainty $\Delta$.\normalsize}
\vspace*{-2pt}
\label{fig:LFT}
\end{center}
\end{figure}
IQC analysis is based on Integral quadratic constraints (IQCs) which, in this context, allow us to describe the uncertainty in the system. Particularly, it is said that a bounded and causal operator $\Delta : \mathbb R^d \rightarrow \mathbb R^d$ satisfies the IQC defined by $\Pi$, i.e., $\Delta \in \IQC(\Pi)$, if 

\vspace*{-15pt}
\small
\begin{align}\label{eq:IQCT}
\int_{0}^{\infty} \begin{bmatrix} v \\ \Delta(v) \end{bmatrix}^{T} \Pi \begin{bmatrix} v \\ \Delta(v) \end{bmatrix} \, dt \geq 0, \quad \forall v \in \mathcal{L}_2^d \ , 
\end{align}
\normalsize
where $\Pi$ is a bounded and self adjoint operator. This constraint can also be rewritten in the frequency domain as 

\vspace*{-16pt}
\small
\begin{align}\label{eq:IQCF}
\int_{-\infty}^{\infty} \begin{bmatrix} \widehat{v}(j\omega) \\ \widehat{\Delta(v)}(j\omega) \end{bmatrix}^{\ast} \Pi(j\omega) \begin{bmatrix} \widehat{v}(j\omega) \\ \widehat{\Delta(v)}(j\omega) \end{bmatrix} \, d\omega \geq 0, 
\end{align} 
\normalsize
where $\hat v$ and $\widehat{\Delta(v)}$ are the Fourier transforms of the signals. Consider the following uncertain system

\vspace*{-16pt}
\small
\begin{equation}\label{eq:UncertainSystem}
\begin{split}
p &= G q\\ q &= \Delta(p),
\end{split}
\end{equation}
\normalsize
where $G \in \mathcal{RH}_{\infty}^{m\times m}$ is the so-called system transfer function matrix and $\Delta : \mathbb R^m \rightarrow \mathbb R^m$ (a bounded and causal operator) is the uncertainty in the system, see Figure~\ref{fig:LFT}. Assume that $\Delta \in \IQC(\Pi)$. Then under certain technical assumptions, see \cite[Thm. 1]{rantzer}, the uncertain system in ~\eqref{eq:UncertainSystem} is robustly stable if 

\vspace*{-16pt}
\small
\begin{align}\label{eq:thmIQC}
\begin{bmatrix} G(j\omega) \\ I \end{bmatrix}^{\ast} \Pi(j\omega) \begin{bmatrix} G(j\omega) \\ I \end{bmatrix} \preceq -\epsilon I,  \hspace{2mm} \forall \ \omega \in [0, \infty].
\end{align}
\normalsize
%
%\begin{IEEEproof}
%See \cite{ulfiqc,rantzer}.
%\end{IEEEproof}
%
In order to establish this condition, we are required to find a multiplier $\Pi$, such that $\Delta \in \IQC(\Pi)$ and such that the semi-infinite LMI in \eqref{eq:thmIQC} is satisfied. Commonly, the condition $\Delta \in \IQC(\Pi)$ imposes structural constraints on $\Pi$. The robust stability analysis problem will then boil down to finding a multiplier, with the required structure, such that the LMI in \eqref{eq:thmIQC} is satisfied for all frequencies. There are mainly two approaches for solving this LMI,  which are based on frequency gridding and the use of KYP lemma. The frequency-gridding-based approach solves the analysis problem approximately, by establishing the feasibility of the LMI in \eqref{eq:thmIQC} for finite number of frequencies. This approach preserves the structure in the LMI which, as we will see in Section \ref{sec:Chordal}, will enable us to solve the analysis problem for interconnected systems, efficiently. Next, we provide a definition for interconnected uncertain systems.
\vspace*{-15pt}
%=======================================================================================
\section{Interconnected Systems: A Definition}
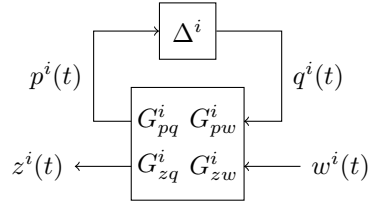
\begin{figure}
\begin{center}
\begin{tikzpicture}[scale=0.5, inner sep=0pt]

 \draw (0,0) node [draw, fill=none, minimum size=1.5cm] {};
 \draw (0,3) node (D) [draw, fill=none, minimum size=0.75cm] {$\Delta^i$};

 \draw[->] (-1.5,0.6) node[right=2pt]{$G^i_{pq}$} -- (-2.5,0.6) -- node[left=4pt]{$p^i(t)$} (-2.5,3) -- (D);
 \draw[->] (D) -- (2.5,3) -- node[right=4pt]{$q^i(t)$} (2.5,0.6) -- (1.5,0.6) node[left=2pt]{$G^i_{pw}$};
 \draw[->] (-1.5,-0.6) node[right=2pt]{$G^i_{zq}$} -- (-3,-0.6) node[left=4pt]{$z^i(t)$};
 \draw[->] (3,-0.6) node[right=4pt]{$w^i(t)$} -- (1.5,-0.6) node[left=2pt]{$G^i_{zw}$};

\end{tikzpicture}  
\vspace*{-7pt}
\caption{\small Individual uncertain subsystems in the network.\normalsize}
%\vspace*{-2pt}
\label{fig:LFTUncertain}
\end{center}
\end{figure}

%=======================================================================================
\vspace*{-10pt}
There are different definitions of interconnected uncertain systems. In this paper, we focus on the following description of interconnected uncertain systems. Let each subsystem in the network be described as
\begin{equation}\label{eq:Subsystems}
\begin{split}
&p^i = G_{pq}^i q^i + G_{pw}^iw^i \\
&z^i = G_{zq}^i q^i + G_{zw}^iw^i\\
&q^i = \Delta^i(p^i),
\end{split}
\end{equation}
where $G_{pq}^i \in \mathcal{RH}_{\infty}^{d_i \times d_i}$, $G_{pw}^i \in \mathcal{RH}_{\infty}^{d_i \times m_i}$, $G_{zq}^i \in \mathcal{RH}_{\infty}^{l_i \times d_i}$, $G_{zw}^i \in \mathcal{RH}_{\infty}^{l_i \times m_i}$, and $\Delta^i:\mathbb{R}^{d_i} \to \mathbb{R}^{d_i}$ is the uncertainty in the $i$th subsystem with $\Delta^i \in \IQC(\Pi^i)$, where $\Pi^i = \begin{bmatrix} \Pi^i_{11} & \Pi^i_{12} \\ \Pi^i_{21} & \Pi^i_{22} \end{bmatrix}$. This is also illustrated in Figure \ref{fig:LFTUncertain}. We define the interconnection among subsystems using the so-called interconnection constraint given by

\vspace*{-13pt}
\small\begin{align}\label{eq:Interconst}
  \underbrace{ \begin{bmatrix}
    w^1\\w^2\\ \vdots\\ w^N
  \end{bmatrix}}_{w} = 
\underbrace{  \begin{bmatrix}
    \Gamma_{11} & \Gamma_{12} & \cdots & \Gamma_{1N} \\
    \Gamma_{21} & \Gamma_{22} & \cdots & \Gamma_{2N} \\
    \vdots & \vdots & \ddots & \vdots \\
    \Gamma_{N1} & \Gamma_{N2} & \cdots & \Gamma_{NN} 
  \end{bmatrix}}_{\Gamma}
  \underbrace{ \begin{bmatrix}
    z^1\\z^2\\ \vdots\\ z^N
  \end{bmatrix}}_{z}.
\end{align}
\normalsize
where $\Gamma$ is the interconnection matrix which has only 0-1 components, and describes how different components of the input-output vectors of different subsystems are connected to one another. Having defined the interconnections, the interconnected uncertain system can be expressed as 
\vspace*{-4pt}
\begin{equation}\label{eq:SysInter}
\begin{split}
p& = G_{pq} q + G_{pw}w \\
z& = G_{zq} q + G_{zw}w\\
q& = \Delta(p)\\
w& = \Gamma z,
\end{split}
\end{equation}
where $G_{\star\bullet} = \diag(G_{\star\bullet}^1, \dots, G_{\star\bullet}^N)$ and $\Delta  = \diag(\Delta^1, \dots, \Delta^N)$. Notice that the system description given in \eqref{eq:SysInter} can readily accommodate uncertainty in the interconnections among the subsystems. This is described in the next section. 
\vspace*{-8pt}
\subsection{Uncertain Interconnections}
\vspace*{-8pt}
Assume that the interconnection between the $i$th subsystem and its neighbors is uncertain, i.e.,
\vspace*{-3pt}
\begin{align*}
w^i = \tilde \Delta^i(\Gamma_i z),
\end{align*}
where $\Gamma_i$ is the $i$th block row of the interconnection matrix $\Gamma$, and $\tilde \Delta^i \in \IQC (\tilde \Pi^i)$ represents the uncertainty in the interconnection. We can then modify the system description of the $i$th subsystem as below

\vspace*{-13pt}
\small
\begin{equation}
\begin{split}
\bar p^i &= \begin{bmatrix} p^i \\ p^i_A \end{bmatrix} = \underbrace{\begin{bmatrix}  G_{pq}^i & G_{pw}^i \\ 0 & 0 \end{bmatrix}}_{\bar G^i_{\bar p \bar q}}\begin{bmatrix} q^i \\ q^i_A \end{bmatrix}+\underbrace{\begin{bmatrix} 0\\I  \end{bmatrix}}_{\bar G^i_{\bar p w}}w^i \\
z^i & = \underbrace{\begin{bmatrix} G_{zq}^i & G_{zw}^i \end{bmatrix}}_{\bar G^i_{z \bar q}}\begin{bmatrix} q^i \\ q^i_A \end{bmatrix}\\
\bar q^i & = \begin{bmatrix} q^i \\ q^i_A \end{bmatrix} = \underbrace{\begin{bmatrix} \Delta^i & 0 \\ 0 &  \tilde \Delta^i \end{bmatrix}}_{\bar \Delta^i}\begin{bmatrix}  p^i \\ p^i_A\end{bmatrix},
\end{split}
\end{equation}
\normalsize
which, if put together, results in the following description of the interconnected uncertain system 
\begin{equation}\label{eq:UncertainInter}
\begin{split}
\bar p &= \bar G_{\bar p \bar q} \bar q + \bar G_{\bar p w} w\\
z & = \bar G_{z \bar q} \bar q\\
\bar q &= \bar \Delta \bar p\\
w & = \Gamma z.
\end{split}
\end{equation}
Similar to \eqref{eq:SysInter}, $\bar p$ and $\bar q$ are the $\bar p^i$s and $\bar q^i$s stacked and the transfer function matrices and $\bar \Delta$ are all block-diagonal. With some re-ordering, the description given in \eqref{eq:UncertainInter} can also be rewritten as

\vspace*{-13pt}
\small
\begin{equation}\label{eq:UncertainInter1}
\begin{split}
\begin{bmatrix}p \\ p_A \end{bmatrix} &= \begin{bmatrix}G_{pq} & G_{pw}\\ 0 & 0  \end{bmatrix}\begin{bmatrix}q \\ q_A \end{bmatrix}+\begin{bmatrix}0 \\ I \end{bmatrix}w\\
z &= \begin{bmatrix} G_{zq} & G_{zw} \end{bmatrix}\begin{bmatrix}q \\ q_A \end{bmatrix}\\
\begin{bmatrix}q \\ q_A \end{bmatrix} &= \begin{bmatrix} \Delta & 0\\0 & \tilde \Delta  \end{bmatrix} \begin{bmatrix}p \\ p_A \end{bmatrix}\\
w & = \Gamma z.
\end{split}
\end{equation}
\normalsize
We will next use the descriptions in \eqref{eq:SysInter} and \eqref{eq:UncertainInter1} to analyze interconnected uncertain systems (with certain and uncertain interconnections) and will study the structure in the resulting LMIs.
\vspace*{-5pt}
%=======================================================================================
\section{Robust Stability Analysis of Interconnected Uncertain Systems}
%=======================================================================================
\vspace*{-5pt}
%
%\subsection{Frequency Gridding Based Approach}
We first start by considering the interconnected system description given in \eqref{eq:SysInter}, i.e., no uncertainty in the interconnections. It is possible to eliminate the interconnection constraint and rewrite the description as
\begin{equation}\label{eq:Lumped}
\begin{split} 
p = \bar G q,\quad q = \Delta(p),
\end{split}
\end{equation}
where $\bar G = G_{pq} + G_{pw}(I - \Gamma G_{zw})^{-1}\Gamma G_{zq}$ and $\Delta \in \IQC(\bar \Pi)$ with 
\begin{align}\label{eq:IQCDiag}
\bar{\Pi} = \begin{bmatrix} \bar{\Pi}_{11} & \bar{\Pi}_{12}\\ \bar{\Pi}_{21} & \bar{\Pi}_{22}\end{bmatrix},
\end{align}
and $\bar \Pi_{ij} = \diag( \Pi_{ij}^1, \dots,  \Pi_{ij}^N)$. We will refer to $\bar G$ as the lumped system transfer function matrix. Note that $I - \Gamma G_{zw}$ must have a bounded inverse for all frequencies in order for the interconnection to be well-posed. This puts the interconnected system in the form in \eqref{eq:UncertainSystem}, and enables us to directly use frequency gridding to analyze the system. We refer to this formulation of the analysis problem as the lumped formulation. However, notice that even if the interconnection among the subsystems is sparse, i.e., $\Gamma$ is sparse, the lumped system transfer matrix will not necessarily be sparse. This is because even though $I - \Gamma G_{zw}$ is sparse, its inverse will generally be dense, which can be seen from the Cayley-Hamilton theorem. As a result, analyzing the large-scale interconnected uncertain system using its lumped description requires solving a large dense LMI, which can be prohibitively costly. In order to alleviate this issue, in \cite{and:13,and+han:12}, we put forth an alternative way of formulating the analysis problem which allowed us to avoid forming the dense lumped system and exploit the inherent sparsity in the problem. Particularly, it was shown in those papers that we can, without additional conservativeness, solve the analysis problem by checking the feasibility of the following LMI

\vspace*{-14pt}
\small
\begin{multline}\label{eq:IQCInterconnected}
\begin{bmatrix} G_{pq} & G_{pw} \\ I & 0 \end{bmatrix}^{\ast}\begin{bmatrix} \bar{\Pi}_{11} & \bar{\Pi}_{12} \\ \bar{\Pi}_{21} & \bar{\Pi}_{22} \end{bmatrix}\begin{bmatrix} G_{pq} & G_{pw} \\ I & 0 \end{bmatrix} -\\
 \begin{bmatrix} -G_{zq}^{\ast}\Gamma^T\\ I
  -G_{zw}^{\ast}\Gamma^T \end{bmatrix}X\begin{bmatrix} -\Gamma
  G_{zq} & I-\Gamma G_{zw} \end{bmatrix} \preceq -\epsilon I.
\end{multline} 
\normalsize
In \eqref{eq:IQCInterconnected} $X$ can be chosen to be a diagonal matrix and hence, this LMI will be sparse in case $\Gamma$ is sufficiently sparse and can then often be solved efficiently using sparse solvers in a centralized manner, \cite{and:13}. We will later in Section \ref{sec:Chordal} show how the sparsity in \eqref{eq:IQCInterconnected}, allow us to decompose the problem. \newline
However, For analyzing interconnected uncertain systems with uncertain interconnections, the lumped formulation of the analysis problem will already be sparse. This is because for the description in \eqref{eq:UncertainInter1}, the lumped system transfer function matrix $\tilde G$ is given by

\vspace*{-14pt}
\small
\begin{align}
\tilde G  =\begin{bmatrix}G_{pq} & G_{pw}\\ 0 & 0  \end{bmatrix} + \begin{bmatrix}0 \\ I \end{bmatrix}\Gamma\begin{bmatrix} G_{zq} & G_{zw} \end{bmatrix} = \begin{bmatrix} G_{pq} & G_{pw}\\ \Gamma G_{zq} & \Gamma G_{zw}   \end{bmatrix},
\end{align}
\normalsize
which is sparse if $\Gamma$ is sparse. As a result, if we directly apply the IQC analysis to this description we will arrive at a sparse LMI. Particularly, assume that $\tilde \Delta = \diag (\tilde \Delta^1, \dots, \tilde \Delta^N) \in \IQC (\tilde \Pi)$, where $\tilde \Pi$ has the same format as in \eqref{eq:IQCDiag}. The resulting LMI from applying the IQC theorem to this system will then take the following form 

\vspace*{-15pt}
\small
\begin{align*}
\begin{bmatrix}  G_{pq} & G_{pw}\\ \Gamma G_{zq} & \Gamma G_{zw} \\ I & 0\\ 0 & I \end{bmatrix}^{\ast} \begin{bmatrix} \bar \Pi_{11} & 0 & \bar \Pi_{12} & 0\\ 0 & \tilde \Pi_{11} & 0 & \tilde \Pi_{12} \\ \bar \Pi_{21} & 0 & \bar \Pi_{22} & 0 \\0 & \tilde \Pi_{21} & 0 & \tilde \Pi_{22}  \end{bmatrix} \begin{bmatrix}  G_{pq} & G_{pw}\\ \Gamma G_{zq} & \Gamma G_{zw} \\ I & 0\\ 0 & I \end{bmatrix} \preceq -\epsilon I,
\end{align*}
\normalsize
for all $\omega \in [0, \infty]$, which can be rewritten as 

\vspace*{-17pt}
\footnotesize
\begin{multline}\label{eq:sparseA}
\begin{bmatrix}G_{pq}^* \\ G_{pw}^* \end{bmatrix}\bar \Pi_{11} \begin{bmatrix} G_{pq} & G_{pw} \end{bmatrix} + \begin{bmatrix}  G_{zq}^*\Gamma^T \\  G_{zw}^*\Gamma^T \end{bmatrix} \tilde \Pi_{11} \begin{bmatrix}\Gamma G_{zq} & \Gamma G_{zw} \end{bmatrix} +\\ \begin{bmatrix}\bar \Pi_{21} & 0 \\0 & \tilde \Pi_{21}\end{bmatrix} \begin{bmatrix} G_{pq} & G_{pw}\\ \Gamma G_{zq} & \Gamma G_{zw} \end{bmatrix}+ \begin{bmatrix} G_{pq} & G_{pw}\\ \Gamma G_{zq} & \Gamma G_{zw} \end{bmatrix}^* \begin{bmatrix}\bar \Pi_{12} & 0 \\0 & \tilde \Pi_{12}\end{bmatrix} +\\ \begin{bmatrix}\bar\Pi_{22} & 0 \\0 & \tilde \Pi_{22} \end{bmatrix} \preceq -\epsilon I. 
\end{multline}
\normalsize
This LMI will also be sparse if $\Gamma$ is sufficiently sparse. 
\begin{rem}
Notice that we can use the system description in \eqref{eq:UncertainInter1} also for certain interconnections, i.e., $\tilde \Delta = I$, and hence \eqref{eq:sparseA} can also be used for analyzing such systems. This actually results in the same LMI as in \eqref{eq:UncertainInter1}. Similarly the formulation in \eqref{eq:IQCInterconnected} can be used to analyze the system in \eqref{eq:UncertainInter1}. However, It then requires solving an LMI with larger dimensions than of \eqref{eq:sparseA}.
\end{rem}
\vspace*{-2pt}
The LMIs in \eqref{eq:IQCInterconnected} and \eqref{eq:sparseA} constitute semidefinite programs (SDPs), \cite{boyd:04}, and next, we will see how the sparsity in these problems can be used to decompose these SDPs and solve them distributedly. 

%%%%%%%%%%%%%%%%%%%%%%%%%%%%%%%%%%%%%%%
\vspace*{-5pt}
\section{Chordal Graphs and Sparsity in SDPs}\label{sec:Chordal}
\vspace*{-8pt}
Chordal sparsity plays a fundamental role in many sparse matrix algorithms, \cite{blp:94}. In Section~\ref{sec:ChordalGraph}, we first review some of the fundamental concepts in graph theory, that will later be used for expressing and exploiting chordal sparsity in SDPs in Section \ref{sec:ChordalSDP}.
%%%%%%%%%%%%%%%%%%%%%%%%%%%%%%%%%%%%%%%%%%%%
\vspace*{-8pt}
\subsection{Chordal Graphs}\label{sec:ChordalGraph}
\vspace*{-8pt}
Let $Q(V,E)$ represent a graph on its set of vertices $V = \{v_1, \dots, v_n\}$ with $E \subseteq V\times V$ denoting its set of edges. Vertices $v_i, v_j \in V$ are adjacent if $(v_i,v_j) \in E$, and we denote the set of adjacent vertices of $v_i$ by $\adj (v_i) = \{ v_j \in V | (v_i, v_j) \in E \}$. The degree of a vertex in a graph is defined as the number of its adjacent vertices, i.e., $\degree(v_i) = |\adj(v_i)|$. The adjacency matrix of a graph $Q(V,E)$ is defined as a $|V| \times |V|$ matrix $A$ where

\vspace*{-15pt}
\small
\begin{align}
A_{ij} = \begin{cases}
1  \quad (i,j) \in E \\
0  \quad (i,j) \notin E
\end{cases}.
\end{align}
\normalsize
A graph is said to be complete if all its vertices are adjacent. An induced graph by $V^\prime \subseteq V$ on $Q(V,E)$, is a graph $Q_I(V^\prime,E^\prime)$ where $E^\prime = E\cap V^\prime \times V^\prime$. A clique $C_i$ of $Q(V,E)$ is a maximal subset of $V$ that induces a complete subgraph on $Q$, i.e., no clique is properly contained in another clique, \cite{blp:94}. Assume that all cycles of length at least four of $Q(V,E)$ have a chord, where a chord is an edge between two non-consecutive vertices in a cycle. The graph is then called chordal \cite[Ch. 4]{gol:04}. Let $C_Q = \{ C_1, \dots, C_l \}$ denote the set of its cliques, where $l$ is the number of cliques of the graph. Then there exists a tree defined on $C_Q$ such that for every $C_i, C_j \in C_Q$ where $i \neq j$, $C_i \cap C_j$ is contained in all the cliques in the path connecting the two cliques in the tree. This property is called the clique intersection property, \cite{blp:94}. Trees with this property are referred to as clique trees. Another fundamental property of chordal graphs is the so-called running intersection property. Define 

\vspace*{-17pt}
\small
\begin{equation}\label{eq:RS}
\begin{split}
& R_1 = C_1, \quad R_i = C_i \setminus \bigcup_{j=1}^{i-1} C_j,\\
& S_1 = \emptyset, \hspace{6mm} S_i = C_i \cap \bigcap_{j=1}^{i-1}C_{j},
\end{split}
\end{equation}
\normalsize
where $R_i$s and $S_i$s are referred to as residuals and separators of the clique tree. Then the running intersection property states that for every chordal graph, there exists an ordering of its cliques $C_Q = \{ C_1, \dots, C_l \}$ where $C_1$ is the root of the clique tree and each parent in the tree has an index smaller than its children such that 

\vspace*{-15pt}
\small
\begin{equation}\label{eq:RIP}
\begin{split}
&R_i = C_i \setminus C_{\parent(i)}\\ 
&S_i = C_i \cap C_{\parent(i)},
\end{split}
\end{equation}
\normalsize
for $i >1$, where $C_{\parent(i)}$ is the parent of $C_i$ in the clique tree, \cite{and:11}. Representing chordal graphs using clique trees with running intersection property proves to be beneficial in handling optimization problems with chordal sparsity, \cite{blp:94}, which is the subject of the upcoming subsection.
%%%%%%%%%%%%%%%%%%%%%%%%%%%%%%%%%%%%%%%%%%%%
\vspace*{-5pt}
\subsection{Chordal sparsity in SDPs}\label{sec:ChordalSDP}
\vspace*{-8pt}

Graphs can be used to describe the sparsity pattern of matrices. Particularly, the sparsity pattern of a matrix $X \in \mathbf S^{n\times n}$ can be represented as a graph $Q(V,E)$ where $V = \{ 1, \dots, n \}$ and $E = \{ (i,j)\in V\times V \ | \ \ X_{ij} \neq 0, \ i \neq j  \}$. A matrix is then said to be chordal if its sparsity graph is chordal. Consider the following standard form SDP

\vspace*{-15pt}
\small
\begin{subequations}\label{eq:SDPDual}
\begin{align}
 \minimize & \quad  b^T y \label{eq:SDPDual1} \\
 \subject & \quad  \sum_{i = 1}^{m} y_iQ^i + W  \preceq 0 \label{eq:SDPDual2}
 \end{align}
\end{subequations}
\normalsize
where $y \in \mathbb R^{m}$ and $W, Q^i \in \mathbf S^{n}$ for $i = 1, \dots, m$. This optimization problem can be decomposed based on its sparsity structure. Assume that the data matrices of the optimization problem in \eqref{eq:SDPDual}, i.e., $W$ and $Q^i$ for $i = 1, \dots, m$, are chordal with corresponding graph $Q(V,E)$, where $\{ C_1, \dots, C_l \}$ is the set of its cliques with ordering that satisfies the running intersection property. The following theorem plays a central role in decomposition of the problem in~\eqref{eq:SDPDual}.

 %The cliques of this graph correspond to the maximal dense principle submatrices of $X$. 
 %
 \begin{thm}(\cite[Thm. 2.3]{ahmr:88}) \label{thm:split}
Let $Q(V,E)$ be a chordal graph with $C_Q = \{ C_1, \dots, C_l\}$ denoting its cliques set with ordering that satisfies the running intersection property. Then for any negative semidefinite $H \in \mathbf S^n$ with corresponding graph $Q$, there exist negative semidefinite matrices $H^i \in \mathbf S^{|C_k|}$ for $i = 1, \dots, l$ such that
\vspace*{-13pt}
\small
 \begin{align}\label{eq:thmchord}
 H = \sum_{i=1}^{l} E_i H^i E_i^T,
 \end{align}
 \normalsize
 where given $C_i = \{ j_1, \dots, j_r \}$ with $j_1 < j_2 < \dots < j_r$, $E_i = \begin{bmatrix} e_{j_1} & \dots & e_{j_r} \end{bmatrix} $ with $e_j$ representing the $j$th column of $I_n$.
 \end{thm}
\vspace*{-8pt}
In other words, a chordal matrix, $H$, can be decomposed into a sum of matrices as in \eqref{eq:thmchord}, such that $H$ is negative semidefinite if and only if all $H^i$ for $i = 1, \dots, l$, are negative semidefinite, \cite[Sec. 5.1]{kim+koj+mev+yam10, Kak:10}. Define $Q(y) = \sum_{i=1}^{m} y_iQ^i$. Then by Theorem \ref{thm:split} and given $Q(y) = \sum_{i=1}^{l}E_i \bar Q^i(y) E_i^T$ and $W = \sum_{i=1}^{l}E_iW^iE_i^T$, we can rewrite the optimization problem in \eqref{eq:SDPDual} as below

\vspace*{-15pt}
\small
\begin{subequations}\label{eq:MSDPDual}
\begin{align}
 \minimize &\quad  \sum_{i=1}^{m} b_iy_i  \\
 \subject & \quad  \bar Q^k(y) + W^k + U^k(d) \preceq 0, \hspace{3mm} k = 1, \dots, l,
 \end{align}
\end{subequations}
\normalsize
where

\vspace*{-15pt}
\small
\begin{multline}
 U^k(d) =\sum_{i,j \in S_k, \ i \geq j}d_{ijk}E_k^TE_{ij}E_k -\\ \sum_{q \in \children(k)} \ \sum_{i,j \in S_q, \ i \geq j}d_{ijq}E_k^TE_{ij}E_k,
\end{multline}
\normalsize
with $\children(i)$ representing the indices of the  children of the $i$th clique in the clique tree, $d$ a $\bar m$-dimensional vector with components defined as $\{ d_{ijk} \ | \ i,j \in S_k, \ i \geq j, \ k = 2, \dots, l \}$ where $ \bar m = \frac{1}{2}\sum_{k = 1}^{l} |S_k|(|S_k| + 1)$, and

\vspace*{-15pt}
\small
\begin{align}
E_{ij} = \begin{cases}
e_ie_i^T  \hspace{13mm} i = j\\
e_ie_j^T + e_je_i^T \quad i \neq j 
\end{cases}.
\end{align}
\normalsize
This approach for reformulating the SDP in \eqref{eq:SDPDual}, is referred to as the range space decomposition. Note that range space decomposition can be used even if the data matrices are not chordal. In that case, one can consider a so-called chordal embedding of the data matrices which can be computed efficiently, \cite{ocf:76}. By using the range space decomposition method, the semidefinite constraints in \eqref{eq:SDPDual2} are equivalently reformulated as several smaller sized semidefinite constraints together with some additional variables. In case the newly generated semidefinite constraints are of low order, which happens if the cliques contain a small number of elements, and the number of added variables are moderate, which happens in case the overlaps between the cliques are small, the newly formulated problem in \eqref{eq:MSDPDual} can be solved with significantly lower effort in comparison to the problem in~\eqref{eq:SDPDual}. The problems in \eqref{eq:IQCInterconnected} and \eqref{eq:sparseA} are semidefinite feasibility problems, i.e., problems that are described as in~\eqref{eq:SDPDual} with $b=0$. Using the range space decomposition, this feasibility problem can then be rewritten in the following format
\begin{align*}
 \find & \quad  v\\
 \subject & \quad v \in \mathcal C_i, \ \ i = 1, \dots, l 
\end{align*}
where $v = (y,d)$ is a $(m + \bar m)$-dimensional vector, 
\begin{align*}
 \mathcal C_i \equiv \{ v  \,|\, f_i(v) \preceq 0\}, \quad i =1,\ldots,l,
\end{align*}
and $f_i(v) = \bar Q^i(y) + W^i + U^i(d)$. Notice that in this case, $f_i$ only depends on a subset of the variables in vector $v$. Particularly, the number of variables in this subset can be at most $(m + \bar m_i)$ where $ \bar m_i = \frac{1}{2}\big(|S_i|(|S_i| + 1) + \sum_{k\in \children (i)} |S_k|(|S_k| + 1)\big)$. Then in case $\bar m_i \ll \bar m$ and if the coupling among the constraints $\mathcal C_i$ is small, the analysis problem can be solved efficiently using distributed computational tools as presented in \cite{kho:13}. Next, we apply the proposed decomposition approach to an interconnection of 500 uncertain subsystems and study its performance.
\vspace*{-6pt}
\section{Numerical Experiment}
\begin{figure}[t]
\begin{center}
\includegraphics[width=6.9cm]{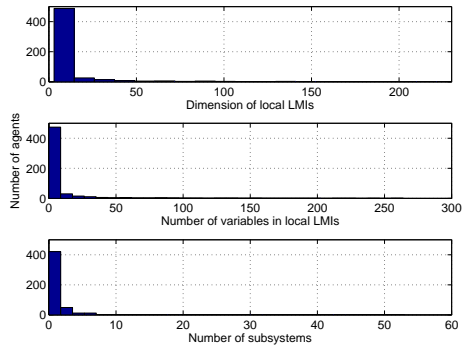}    % The printed column width is 8.4 cm.
\vspace*{-18pt}
\caption{\small The histogram of dimension and number of variables of each of the constraints (top and middle plots). The plot on the bottom, shows different agents required information from how many subsystems.\normalsize }
\vspace*{-3pt}
\label{fig:Dim}
\end{center}
\end{figure}

\vspace*{-8pt}
In this section, we analyze the robustness of an interconnection of 500 uncertain subsystems. In order to define the interconnection among the subsystems and each subsystems' transfer function matrix, we have used the same procedure that was described in \cite{and:13}. The considered adjacency matrix of the network, that describe the interconnection among the subsystems, is such that only 6 of the subsystems are connected to more than 11 other subsystems, which implies sparsity in the interconnections among the subsystems.  The uncertainty block in each subsystem has dimension 1 and corresponds to parametric uncertainty in the subsystem and is assumed to be an unknown gain normalized within $[-1,1]$. Analyzing this interconnected system using the lumped formulation requires solving a fully dense LMI with dimension 500 and 500 variable. Performing this analysis using the sparse formulation in \eqref{eq:IQCInterconnected}, however, requires solving a very sparse LMI of dimension 1498 with 1498 variables, where $X$ is chosen to be a diagonal matrix variable. The chordal embedding of this LMI has 579 cliques and hence it can be decomposed as in \eqref{eq:MSDPDual}, which results in a problem with 579 semidefinite constraints. This problem includes 9894 variables. The histogram of dimension and number of variables of each of these constraints are presented in Figure \ref{fig:Dim}. As can be seen from the figure, most of the constraints, in fact 94\% of them, have dimension and number of variables less than 50. Also the largest semidefinite constraint has dimension 210 and has 170 variables, which is far less than 9894. The coupling among the constraints is also small, where 95 \% of the constraints have couplings that include less than 24 variables. The largest coupling includes 92 variables. This implies that the problem can be solved efficiently using the algorithms presented in \cite{kho:13} via 579 computing agents, and the hardest problem that has to be solved by these agents requires solving a semidefinite program of dimension 210 and with 170 variables. For this numerical experiment we have utilized Algorithm 9 in \cite{kho:13}, which converged to a feasible solution within 14 iterations. Notice that the computational cost of solving the lumped formulation is not prohibitive and can be solved within 2760 seconds in a centralized manner. This is in comparison to 1623 seconds for the distributed solution with 579 parallel processors. With the distributed solution, however, we have also been able to provide partial local privacy within the network, in the sense that each agent, in order to carry out its local computations, requires information about the model of only a few of the subsystems. This can be seen from the bottom plot of Figure \ref{fig:Dim}, where it shows that 96~\% of the agents require information from less than 10 subsystems, and in the worst case there is 1 agent that requires information from 52 subsystems.   
\vspace*{-8pt}
\section{Conclusions}
\vspace*{-9pt}
In this paper, we considered robust stability analysis of large-scale interconnected uncertain systems, with certain and uncertain interconnections, and proposed a distributed scheme for solving the analysis problem. This was achieved by observing that the analysis of interconnected uncertain systems can be conducted by solving sparse LMIs, sparsity pattern of which depends on the interconnection structure among the subsystems. In case the interconnection matrix is sufficiently sparse, this will then enable us to decompose these LMIs and solve them distributedly using distributed computational algorithms. As future research, we intend to study other sparsity-based decomposition techniques and investigate how the decomposition can be done more efficiently to reduce the computational cost of solving the problem distributedly. Also we will apply the proposed decomposition scheme to larger problems, that cannot be solved in centralized manner, to highlight the method's computational benefits.
\vspace*{-17pt}
\small
\bibliography{IEEETrans}

\begin{thebibliography}{28}
\providecommand{\natexlab}[1]{#1}
\providecommand{\url}[1]{\texttt{#1}}
\expandafter\ifx\csname urlstyle\endcsname\relax
  \providecommand{\doi}[1]{doi: #1}\else
  \providecommand{\doi}{doi: \begingroup \urlstyle{rm}\Url}\fi

\bibitem[Agler et~al.(1988)Agler, Helton, McCullough, and Rodman]{ahmr:88}
J.~Agler, W.~Helton, S.~McCullough, and L.~Rodman.
\newblock Positive semidefinite matrices with a given sparsity pattern.
\newblock \emph{Linear Algebra and its Applications}, 107:\penalty0 101--149,
  1988.

\bibitem[Andersen(2011)]{and:11}
M.~S. Andersen.
\newblock \emph{Chordal Sparsity in Interior-Point Methods for Conic
  Optimization}.
\newblock PhD dissertation, university of {California, Los Angeles}, 2011.

\bibitem[Andersen et~al.(2010)Andersen, Dahl, and Vandenberghe]{and+dah+van10}
M.~S. Andersen, J.~Dahl, and L.~Vandenberghe.
\newblock Implementation of nonsymmetric interior-point methods for linear
  optimization over sparse matrix cones.
\newblock \emph{Mathematical Programming Computation}, pages 1--35, 2010.

\bibitem[Andersen et~al.(2012)Andersen, Hansson, Pakazad, and
  Rantzer]{and+han:12}
M.~S. Andersen, A.~Hansson, S.~Khoshfetrat Pakazad, and A.~Rantzer.
\newblock Distributed robust stability analysis of interconnected uncertain
  systems.
\newblock In \emph{{Proceedings of the 51st {IEEE} Conference on Decision and
  Control}}, 2012.

\bibitem[{Andersen} et~al.(2013){Andersen}, {Khoshfetrat Pakazad}, {Hansson},
  and {Rantzer}]{and:13}
M.~S. {Andersen}, S.~{Khoshfetrat Pakazad}, A.~{Hansson}, and A.~{Rantzer}.
\newblock {Robust Stability Analysis of Sparsely Interconnected Uncertain
  Systems}.
\newblock \emph{Provisionally accepted for publication in IEEE Transactions on
  Automatic Control}, 2013.
\newblock URL \url{http://arxiv.org/abs/1311.2602}.

\bibitem[Blair and Peyton(1994)]{blp:94}
J.~R.~S. Blair and B.~W. Peyton.
\newblock An introduction to chordal graphs and clique trees.
\newblock In J.~A. George, J.~R. Gilbert, and J.~W-H. Liu, editors, \emph{Graph
  Theory and Sparse Matrix Computations}, volume~56, pages 1--27.
  Springer-Verlag, 1994.

\bibitem[Boyd and Vandenberghe(2004)]{boyd:04}
S.~Boyd and L.~Vandenberghe.
\newblock \emph{Convex Optimization}.
\newblock Cambridge University Press, 2004.

\bibitem[Fang and Antsaklis(2008)]{P5}
H.~Fang and P.~J. Antsaklis.
\newblock Distributed control with integral quadratic constraints.
\newblock In \emph{Proceedings of the 17th {IFAC} World Congress}, volume~17,
  2008.

\bibitem[Golumbic(2004)]{gol:04}
M.~C. Golumbic.
\newblock \emph{Algorithmic Graph Theory and Perfect Graphs}.
\newblock Elsevier, 2nd edition, 2004.

\bibitem[Hansson and Vandenberghe(2000)]{Anders}
A.~Hansson and L.~Vandenberghe.
\newblock Efficient solution of linear matrix inequalities for integral
  quadratic constraints.
\newblock In \emph{Proceedings of the 39th {IEEE} Conference on Decision and
  Control}, volume~5, pages 5033--5034, 2000.

\bibitem[J{\"{o}}nsson(2001)]{ulfiqc}
U.~J{\"{o}}nsson.
\newblock Lecture notes on integral quadratic constraints, May 2001.

\bibitem[J{\"{o}}nsson and Kao(2010)]{Ulf}
U.~J{\"{o}}nsson and C.~Kao.
\newblock A scalable robust stability criterion for systems with heterogeneous
  {LTI} components.
\newblock \emph{{IEEE} Trans. Autom. Control}, 55\penalty0 (10):\penalty0
  2219--2234, October 2010.

\bibitem[J{\"{o}}nsson et~al.(2007)J{\"{o}}nsson, Kao, and Fujioka]{UlfLetter}
U.~J{\"{o}}nsson, C.~Kao, and H.~Fujioka.
\newblock A {P}opov criterion for networked systems.
\newblock \emph{Systems \& Control Letters}, 56\penalty0 (9--10):\penalty0
  603--610, 2007.

\bibitem[Kakimura(2010)]{Kak:10}
N.~Kakimura.
\newblock A direct proof for the matrix decomposition of chordal-structured
  positive semidefinite matrices.
\newblock \emph{Linear Algebra and its Applications}, 433\penalty0
  (4):\penalty0 819--823, 2010.

\bibitem[Kao and Megretski(2003)]{Kao1}
C.~Kao and A.~Megretski.
\newblock A new barrier function for {IQC} optimization problems.
\newblock In \emph{Proceedings of the American Control Conference}, volume~5,
  pages 4281--4286, June 2003.

\bibitem[Kao et~al.(2004)Kao, Megretski, and J{\"o}sson]{Kao2}
C.~Kao, A.~Megretski, and U.~J{\"o}sson.
\newblock Specialized fast algorithms for {IQC} feasibility and optimization
  problems.
\newblock \emph{Automatica}, 40\penalty0 (2):\penalty0 239--252, 2004.

\bibitem[Kao et~al.(2009)Kao, J{\"o}sson, and Fujioka]{Kao}
C.~Kao, U.~J{\"o}sson, and H.~Fujioka.
\newblock Characterization of robust stability of a class of interconnected
  systems.
\newblock \emph{Automatica}, 45\penalty0 (1):\penalty0 217--224, 2009.

\bibitem[{Khoshfetrat Pakazad} et~al.(2013){Khoshfetrat Pakazad}, {Andersen},
  and {Hansson}]{kho:13}
S.~{Khoshfetrat Pakazad}, M.~S. {Andersen}, and A.~{Hansson}.
\newblock {Distributed Solutions for Loosely Coupled Feasibility Problems Using
  Proximal Splitting Methods}.
\newblock \emph{Provisionally accepted for publication in Optimization Methods
  and Software}, 2013.
\newblock URL \url{http://arxiv.org/abs/1306.6807}.

\bibitem[Kim and Braatz(2012)]{kim:12}
K.~K. Kim and R.~D. Braatz.
\newblock On the robustness of interconnected or networked uncertain linear
  multi-agent systems.
\newblock \emph{20th International Symposium on Mathematical Theory of Networks
  and Systems}, 2012.

\bibitem[Kim et~al.(2010)Kim, Kojima, Mevissen, and
  Yamashita]{kim+koj+mev+yam10}
S.~Kim, M.~Kojima, M.~Mevissen, and M.~Yamashita.
\newblock Exploiting sparsity in linear and nonlinear matrix inequalities via
  positive semidefinite matrix completion.
\newblock \emph{Mathematical Programming}, pages 1--36, 2010.

\bibitem[Langbort et~al.(2004)Langbort, Chandra, and D'Andrea]{Langbort}
C.~Langbort, R.~S. Chandra, and R.~D'Andrea.
\newblock Distributed control design for systems interconnected over an
  arbitrary graph.
\newblock \emph{{IEEE} Trans. Autom. Control}, 49\penalty0 (9):\penalty0
  1502--1519, September 2004.

\bibitem[Lestas and Vinnicombe(2006)]{Vinnicombe}
I.~Lestas and G.~Vinnicombe.
\newblock Scalable decentralized robust stability certificates for networks of
  interconnected heterogeneous dynamical systems.
\newblock \emph{{IEEE} Trans. Autom. Control}, 51\penalty0 (10):\penalty0
  1613--1625, October 2006.

\bibitem[Megretski and Rantzer(1997)]{rantzer}
A.~Megretski and A.~Rantzer.
\newblock System analysis via integral quadratic constraints.
\newblock \emph{{IEEE} Trans. Autom. Control}, 42\penalty0 (6):\penalty0
  819--830, June 1997.

\bibitem[Ohtsuki et~al.(1976)Ohtsuki, Cheung, and Fujisawa]{ocf:76}
T.~Ohtsuki, L.~K. Cheung, and T.~Fujisawa.
\newblock Minimal triangulation of a graph and optimal pivoting order in a
  sparse matrix.
\newblock \emph{Journal of Mathematical Analysis and Applications}, 54\penalty0
  (3):\penalty0 622--633, 1976.

\bibitem[Parrilo(2001)]{Parrilo}
P.~A. Parrilo.
\newblock Outer approximation algorithms for {KYP}-based {LMI}s.
\newblock In \emph{Proceedings of the American Control Conference}, volume~4,
  pages 3025--3028, 2001.

\bibitem[Vandenberghe et~al.(2005)Vandenberghe, Balakrishnan, Wallin, Hansson,
  and Roh]{ChapterBook1}
L.~Vandenberghe, V.~R. Balakrishnan, R.~Wallin, A.~Hansson, and T.~Roh.
\newblock Interior-point algorithms for semidefinite programming problems
  derived from the {KYP} lemma.
\newblock In D.~Henrion and A.~Garulli, editors, \emph{Positive polynomials in
  control}, volume 312, pages 195--238. Springer, Feb 2005.

\bibitem[Wallin et~al.(2009)Wallin, Hansson, and Johansson]{WHJ:09}
R.~Wallin, A.~Hansson, and J.~H. Johansson.
\newblock A structure exploiting preprocessor for semidefinite programs derived
  from the {Kalman-Yakubovich-Popov} lemma.
\newblock \emph{{IEEE} Trans. Autom. Control}, 54\penalty0 (4):\penalty0
  697--704, April 2009.

\bibitem[Zhou et~al.(1997)Zhou, Doyle, and Glover]{robustandoptimal}
K.~Zhou, J.~C. Doyle, and K.~Glover.
\newblock \emph{Robust and Optimal Control}.
\newblock Prentice Hall, 1997.

\end{thebibliography}

%
%\appendix
%\section{A summary of Latin grammar}    % Each appendix must have a short title.
%\section{Some Latin vocabulary}         % Sections and subsections are supported  
                                        % in the appendices.
\end{document}